\documentclass[pdftex,a4paper,12pt]{amsart}

\usepackage{geometry}
\usepackage[cp1250]{inputenc}

\usepackage{alltt}
\usepackage{euler}
\usepackage{amsfonts}
\usepackage{amscd}
\usepackage{graphicx}
\usepackage{lpic}
\usepackage{longtable}
\usepackage{url}
\usepackage{import}
\usepackage{amssymb, amsthm, amsmath}
\usepackage[sc]{mathpazo}
\usepackage{enumitem, textcomp, setspace}

\usepackage{tikz}
\usetikzlibrary{intersections}
\usepackage{pgfplots}
\usepgfplotslibrary{external}

\allowdisplaybreaks

\numberwithin{equation}{section}

\theoremstyle{plain}
\newtheorem{theorem}{Theorem}[section]

\newtheorem{lemma}[theorem]{Lemma}
\newtheorem{corollary}[theorem]{Corollary}

\newtheorem{conjecture}[theorem]{Conjecture}

\theoremstyle{definition}

\newtheorem{remark}[theorem]{Remark}

\usepackage[colorlinks]{hyperref}

\usepackage[all]{hypcap}

\definecolor{internalLink}{rgb}{0,0,0.5}
\definecolor{citeLink}{rgb}{0,0.5,0}
\definecolor{urlLink}{rgb}{0,0.5,0.5}

\hypersetup{
	linkcolor=internalLink,
	citecolor=citeLink,
	urlcolor=urlLink
}

\title[Tabulation of knots up to five triple-crossings and oriented moves]{Tabulation of knots up to five triple-crossings\\ and moves between oriented diagrams}

\author{Micha{\l} Jab{\l}onowski}

\address{Institute of Mathematics, Faculty of Mathematics, Physics and Informatics,\newline University of Gda\'nsk, 80-308 Gda\'nsk, Poland}



\email{michal.jablonowski@gmail.com}

\date{\today}

\begin{document}

\maketitle

\begin{abstract}

We enumerate and show tables of minimal diagrams for all prime knots up to the triple-crossing number equal to five. We derive a minimal generating set of oriented moves connecting triple-crossing diagrams of the same oriented knot. We also present a conjecture about a strict lower bound of the triple-crossing number of a knot related to the breadth of its Alexander polynomial.
\end{abstract}

\section{Introduction}

\noindent
It is known since at least V.F.R. Jones observation in 1999 (in his \emph{Planar Algebras, I} c.f.\cite{Ada21}) that any knot and every link has a diagram where, at each of its multiple-points in the plane, exactly three strands are allowed to cross pairwise transversely. For the very recent survey on this topic see \cite{Ada21}. Such triple-point diagrams and moves on them have been studied in several recent papers, such as \cite{Ada13, ACFIPVWZ14, AHP19, JabTro20, Jab21}.
\par
The triple-crossing number of a knot $K$, denoted here by $c_3(K)$, is defined in analogy to the classical (double-crossing) number, as the least number of triple-crossings for any triple-crossing diagram of $K$. There are lower bounds for the triple-crossing number, in terms of double-crossing number $c_3(K)\geq \frac{1}{3}c_2(K)$, and if $K$ is alternating then $c_3(K)\geq \frac{1}{2}c_2(K)$ (see \cite{Ada13}).
\par
In \cite{Jab21} the author prove the following bound of the triple-crossing number $c_3$ by the canonical genus $g_c$. Let $K$ be a knot. Then $c_3(K)\geq 2\cdot g_c(K).$ It follows from this bound that the triple-crossing number is greater or equal to the breadth of the Alexander polynomial $\Delta$, since it is known that $2\cdot g_c(K)\geq \text{breadth}(\Delta(K))$. We propose a conjecture based on our extensive experiments.

\begin{conjecture}\label{t1}
	Let $K$ be a knot, such that $\Delta(K)$ is not monic. Then $$c_3(K)>\text{breadth}(\Delta(K)).$$
\end{conjecture}

A polynomial is called monic if the coefficient of the highest order term are equal to $\pm 1$. If true, the conjecture immediately gives a sharp enough bound to obtain the exact (unknown) value of the triple-crossing number of many knots (from known upper bounds on the triple-crossing number), such as (giving only for knots with $c_2\leq 13$):\\ \ \\ $9_{3}$, $9_{6}$, $9_{9}$, $9_{16}$, $K11a_{234}$, $K11a_{240}$, $K11a_{263}$, $K11a_{334}$, $K11a_{338}$, $K11a_{355}$, $K11a_{364}$, $K13a_{3092}$, $K13a_{3110}$, $K13a_{3132}$, $K13a_{3377}$, $K13a_{3380}$, $K13a_{4547}$, $K13a_{4558}$, $K13a_{4739}$, $K13a_{4822}$, $K13a_{4828}$, $K13a_{4862}$, $K13a_{4874}$.
\ \\
In Section\;\ref{s2} of this paper, we also derive a minimal generating set of oriented moves connecting triple-crossing diagrams of the same oriented knot. Later, in Section\;\ref{s3} we enumerate and show tables of minimal diagrams for all prime knots up to the triple-crossing number equal to five. We use the knot names used in the \cite{katlas} database package, and for $11$--$14$-crossing knots we use the Hoste-Thistlethwaite database.

\section{Definitions}\label{s1}

The \emph{projection} of a knot or a link $K\subset \mathbb{R}^3$ is its image under the standard projection $\pi:\mathbb{R}^3\to\mathbb{R}^2$ (or into a $2$-sphere) such that it has only a finite number of self-intersections, called \emph{multiple points}, and in each multiple-points each pair of its strands are transverse. 
\par
If each multiple-points of a projection has multiplicity three then we call this projection a \emph{triple-crossing projection}. The \emph{triple-crossing} is a three-strand crossing with the strand labeled $T, M, B$, for top, middle and bottom.
\par
\begin{figure}[h!t]
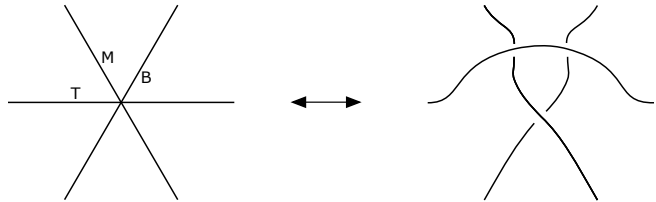

		\centering		
		\begin{lpic}[]{./M01(8.5cm)}
		
		\end{lpic}
		\caption{A deconstruction/construction of a triple-crossing.}
		\label{r01}
	\end{figure}

The \emph{triple-crossing diagram} is a triple-crossing projection such that each of its triple points is a triple-crossing, such that $\pi^{-1}$ of the strand labeled $T$ (in the neighborhood of that triple point) is on the top of the strand corresponding to the strand labeled $M$, and the latter strand is on the top of the strand corresponding to the strand labeled $B$ (see Figure\;\ref{r01}).

The \emph{triple-crossing number} of a knot or link $K$, denoted $c_3(K)$, is the least number of triple-crossings for any triple-crossing diagram of $K$. The classical double-crossing number invariant we will denote by $c_2$. The \emph{minimal triple-crossing diagram} of a knot $K$ is a triple-crossing diagram of $K$ that has exactly $c_3(K)$ triple-crossings.
\par
A \emph{natural orientation} (see an equivalent definition in \cite{AHP19}) on a triple-crossing diagram is an orientation of each component of that link, such that in each crossing the strands are oriented in-out-in-out-in-out, as we encircle the crossing. We begin with an interesting notice.

\begin{lemma}[\cite{AHP19}]\label{l1}
	
	Every orientation of the triple-crossing diagram obtained from an oriented knot is the natural orientation.
\end{lemma}

\section{Oriented moves}\label{s2}

	\begin{figure}[h!t]
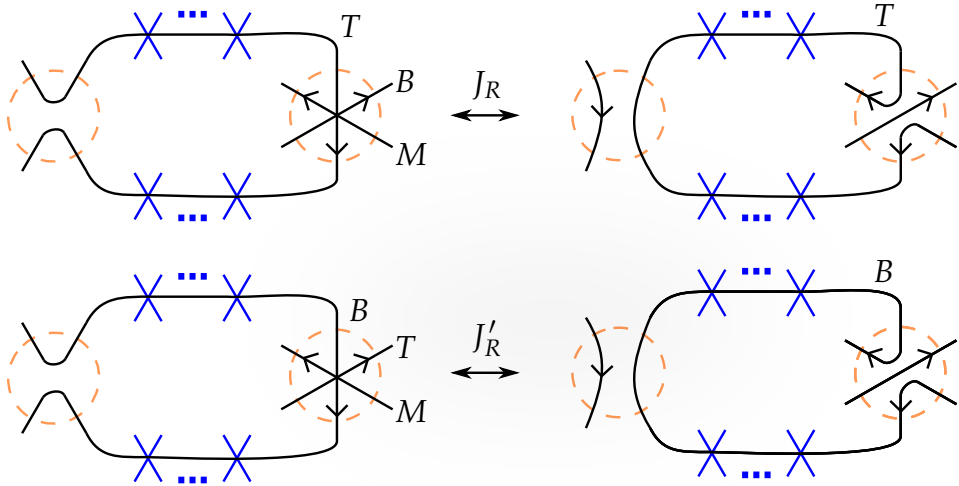

	\centering		
	\begin{lpic}[t(0.2cm)]{./M06ori(12.5cm)}
		\lbl[b]{52,46;$J_R$}
		\lbl[b]{52,18;$J_R'$}
		\lbl[l]{36,54;$T$}
		\lbl[b]{95,54;$T$}
		\lbl[l]{37,23;$B$}
		\lbl[b]{95,26;$B$}
		\lbl[l]{42,12;$M$}
		\lbl[l]{42,40;$M$}
		\lbl[l]{42,48;$B$}
		\lbl[l]{42,19;$T$}
	\end{lpic}
	\caption{A minimal generating set of oriented moves on triple-crossing diagrams of a knot.}
	\label{r06ori}
\end{figure}

\begin{theorem}\label{tw1}
	Two oriented triple-crossing diagrams of knots are related by a sequence of oriented $J_R$ and $J_R'$ moves (see Figure\;\ref{r06ori}) and a spherical isotopy, if and only if they define the same knot type.
\end{theorem}

In the $J_R$ move, there can be finitely many triple-crossings (colored here in blue) such that for every triple-crossing the arc, that is nearest the letter $T$ lies always on top. In the $J_R'$ move that is the mirror move to $J_R$, there can be finitely many triple-crossings (colored here in blue) such that for every triple-crossing the arc, that is nearest the letter $B$ lies always on bottom. The orientations of the blue arcs are determined by the natural orientation property.

\begin{figure}[h!t]
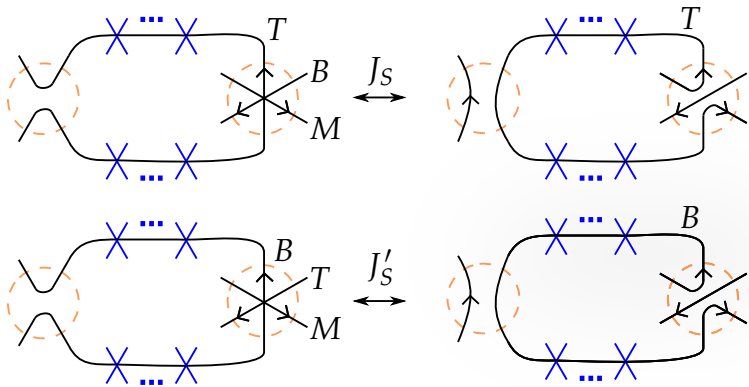

	\centering		
	\begin{lpic}[t(0.2cm)]{./M06oriB(12cm)}
		\lbl[b]{52,46;$J_S$}
		\lbl[b]{52,18;$J_S'$}
		\lbl[l]{36,54;$T$}
		\lbl[b]{95,54;$T$}
		\lbl[l]{37,23;$B$}
		\lbl[b]{95,26;$B$}
		\lbl[l]{42,12;$M$}
		\lbl[l]{42,40;$M$}
		\lbl[l]{42,48;$B$}
		\lbl[l]{42,19;$T$}
	\end{lpic}
	\caption{The other pair of oriented moves.}
	\label{r06oriB}
\end{figure}

\begin{proof}

	We have the minimal set of unoriented moves $J_R$ and $J_R'$ between unoriented knots, defined by the author in \cite{JabTro20} that are identical as our moves in Figure\;\ref{r06ori} but without decorating arrows (so we leave the names unchanged). From Lemma\;\ref{l1} we see that specifying orientation on one strand in any triple-crossing the other strands in that crossing must have determined orientation. Therefore, because in each local diagram of unoriented moves $J_R$ and $J_R'$ there is a strand passing through all other triple-crossings, we have up to four generating moves $J_R$, $J_R'$, $J_S$ and $J_S'$ for oriented diagrams (for the latter pair see Figure\;\ref{r06oriB}). But the move $J_S$ can be generated from $J_R$ by a spherical isotopy.
	First choose any non-outer region adjacent region to the triple-crossings marked $T, M , B$ and on the left to the triple-crossing then by a spherical isotopy make the region to be the outer (unbounded) region for the knot diagram. Then rotate the diagram by $180$ degrees. The same goes with the pair $J_S'$ and $J_R'$.

\end{proof}

\begin{corollary}
	The set $\{J_R, J_R'\}$ (presented in Figure\;\ref{r06ori}) is a minimal generating set of oriented moves connecting triple-crossing diagrams of the same oriented knot.
\end{corollary}

\section{Knot Tabulation}\label{s3}

After several months of computer computations, following the method described by the author in \cite{JabTro20} and implemented in \texttt{Wolfram Engine 12}, we generate the table of knots with the triple-crossing number equal to five.
\par
\begin{figure}[h!t]
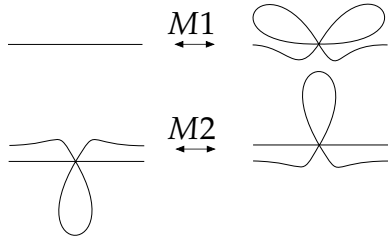

	\centering		
	\begin{lpic}[]{./M02(5cm)}
				\lbl[b]{57,63;$M1$}
				\lbl[b]{57,30;$M2$}
	\end{lpic}
	\caption{Moves $M1$ and $M2$.}
	\label{r02}
\end{figure}

First, we enumerate all prime, connected, triple-crossing projections, up to spherical isotopy, up to mirror image, and up to moves $M1$, $M2$ (see Figure\;\ref{r02}), with five triple points, with the result of $116$ projections (see Table\;\ref{tab1}).

To identify types of knots, we use the classical polynomial invariants, the Jones polynomial (see Remark\;\ref{re1}), and later two-variable polynomials where they are needed.

\begin{remark}\label{re1}
	By Lemma\;\ref{l1} and the Kauffman bracket relations form \cite{Ada13}, the Jones polynomial $V$ for an oriented knot can be calculated from a triple-crossing diagram (regardless of the orientation) by resolving each triple-crossing by the following relations in Table\;\ref{tabJ}.
	
	\begin{table}
		\caption{Relations for the Jones polynomial.}
		\label{tabJ}
		\begin{tabular}{ll}
			
			$V\bigg( \raisebox{-15pt}{\includegraphics[scale=0.3]{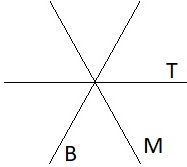}} \bigg)=$&$-t^{\frac{3}{2}}V\bigg( \raisebox{-15pt}{\includegraphics[scale=0.3]{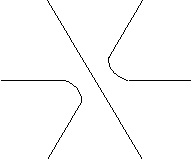}} \bigg)-tV\bigg( \raisebox{-15pt}{\includegraphics[scale=0.3]{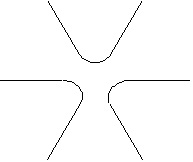}} \bigg)+$\\&\\&$-tV\bigg( \raisebox{-15pt}{\includegraphics[scale=0.3]{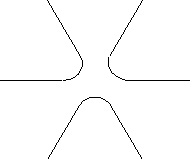}} \bigg)-t^{\frac{1}{2}}V\bigg( \raisebox{-15pt}{\includegraphics[scale=0.3]{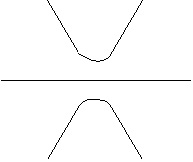}} \bigg)-t^{\frac{1}{2}}V\bigg( \raisebox{-15pt}{\includegraphics[scale=0.3]{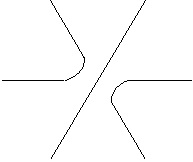}} \bigg)$\\ & \\
			
			$V\bigg( \raisebox{-15pt}{\includegraphics[scale=0.3]{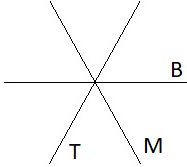}} \bigg)=$&$-t^{-\frac{3}{2}}V\bigg( \raisebox{-15pt}{\includegraphics[scale=0.3]{./sk2}} \bigg)-t^{-1}V\bigg( \raisebox{-15pt}{\includegraphics[scale=0.3]{./sk3}} \bigg)+$\\&\\&$-t^{-1}V\bigg( \raisebox{-15pt}{\includegraphics[scale=0.3]{./sk4}} \bigg)-t^{-\frac{1}{2}}V\bigg( \raisebox{-15pt}{\includegraphics[scale=0.3]{./sk5}} \bigg)-t^{-\frac{1}{2}}V\bigg( \raisebox{-15pt}{\includegraphics[scale=0.3]{./sk6}} \bigg)$\\ & \\
			
			$V(\underbrace{O\sqcup \ldots \sqcup O}_{c})$&$=\left(-t^{\frac{1}{2}}-t^{-\frac{1}{2}}\right)^{c-1}$
			
		\end{tabular}
		
	\end{table}	
\end{remark}

The number of knots with a specific triple-crossing number is presented in Table\;\ref{tab1}. Diagrams of the knots, generated by a new algorithm, are presented in Table\;\ref{tab2}, where the labels are of the form $tk_n$ (for the triple-crossing number equal to $k$, and $n$ the consecutive index in that family). The labels in brackets are the classical Alexander-Briggs-Rolfsen notation of a knot (up to mirror image). The convention here is that the green strand near a triple-crossing is the bottom strand, and the red one is the upper strand (the loops are always black for simplicity).

\begin{table}[ht]
	\caption{Enumeration of knots and projections.}
	\label{tab1}
\begin{tabular}{r||c|c}
	crossings	& number of projections& number of knots\\
		\hline\hline
		$2$&$1$&$2$\\
		$3$&$2$&$2$\\
		$4$&$15$&$24$\\
		$5$& $\bf 116$ &$\bf 118$
	\end{tabular}
\end{table}

\par
An \texttt{sPD} code for a given knot can be found below the \texttt{TikZ} code of the corresponding diagram, in the \texttt{LaTeX} source file of this article’s \texttt{arXiv} preprint version. In that archive there is also a text file of \texttt{sPD} codes of all the mentioned triple-crossing projections.

\begin{longtable}{@{}cccc@{}}
	\caption{Knots with the triple-crossing number $\leqslant 5$.}\\
\end{longtable}%
	
\end{document}